\documentclass[12pt]{article}
\usepackage{amsmath,amssymb,amsthm}

\newtheorem{theorem}{Theorem}

\newtheorem{proposition}{Proposition}

\usepackage{geometry}

\usepackage[shortlabels, inline]{enumitem}

\usepackage{bm}

\usepackage{graphicx}
\usepackage{float}

\allowdisplaybreaks[2]

\usepackage{mathtools}
\mathtoolsset{showonlyrefs=true}

\usepackage{hyperref}

\newcommand{\sankaku}[1]{\textup{$\bigtriangleup$#1}}
\newcommand{\kaku}[1]{\textup{$\angle$#1}}
\newcommand{\ten}[1]{\textup{#1}}

\begin{document}

\begin{center}
\Large
    Appearance of similar triangles by certain operations on triangles

\normalsize

\vspace{15mm}

Hiroki Naka

Nada High School

8-5-1 Uozakikita-cho, Higashinada-ku, Kobe, Hyogo 658--0082, Japan

e111115@stg.nada.ac.jp

\vspace{8mm}

Takahiko Fujita

Department of Science and Engineering, Chuo University,
 
1-13-27, Kasuga, Bunkyo, Tokyo 112-8551, Japan

rankstatistics@gmail.com

\vspace{8mm}

Naohiro Yoshida

Department of Economics, Keiai University,

1-5-21 Anagawa, Inage-ku,
Chiba-shi, Chiba 263-8588, Japan

n-yoshida@u-keiai.ac.jp
    
\end{center}



\begin{abstract}
    In this paper, a theorem about similar triangles is proved. It shows that two small and four large triangles similar to the original triangle can appear if we choose well among several intersections of the perpendicular bisectors of the sides with perpendicular lines of sides passing through the vertices of the triangle.
\end{abstract}

\section{Introduction}

The theorem presented in this paper states that two small and four large triangles similar to the original triangle can appear if we choose well among several intersections of the perpendicular bisectors of the sides with perpendicular lines of sides passing through the vertices of the triangle.
One of the authors, H. Naka, discovered this theorem by chance several years ago in \cite{Naka} while playing with the drawing app named ``GeoGebra", so we will call this triangle \emph{Naka's triangle}. 
In \cite{Naka2022}, he was then able to find a proof of the theorem by elementary geometry and also found the relation with the Miquel points \cite{Miquel}.
In the process, he got acquainted with the second author, T. Fujita, and the third author, N. Yoshida, and started to collaborate.
We were able to give a proof by complex numbers, and obtained complex identities appearing in the calculations, which had not been known before.

The rest of the paper is organized as follows.
Section 2 states our main result.
In Section 3, a proof by complex numbers is given.
Section 4 provides a proof by elementary geometry.
Section 5 is devoted to further extensions.

\section{Main result}

\begin{figure}[t]
    \centering
    \includegraphics[width=0.7\linewidth]{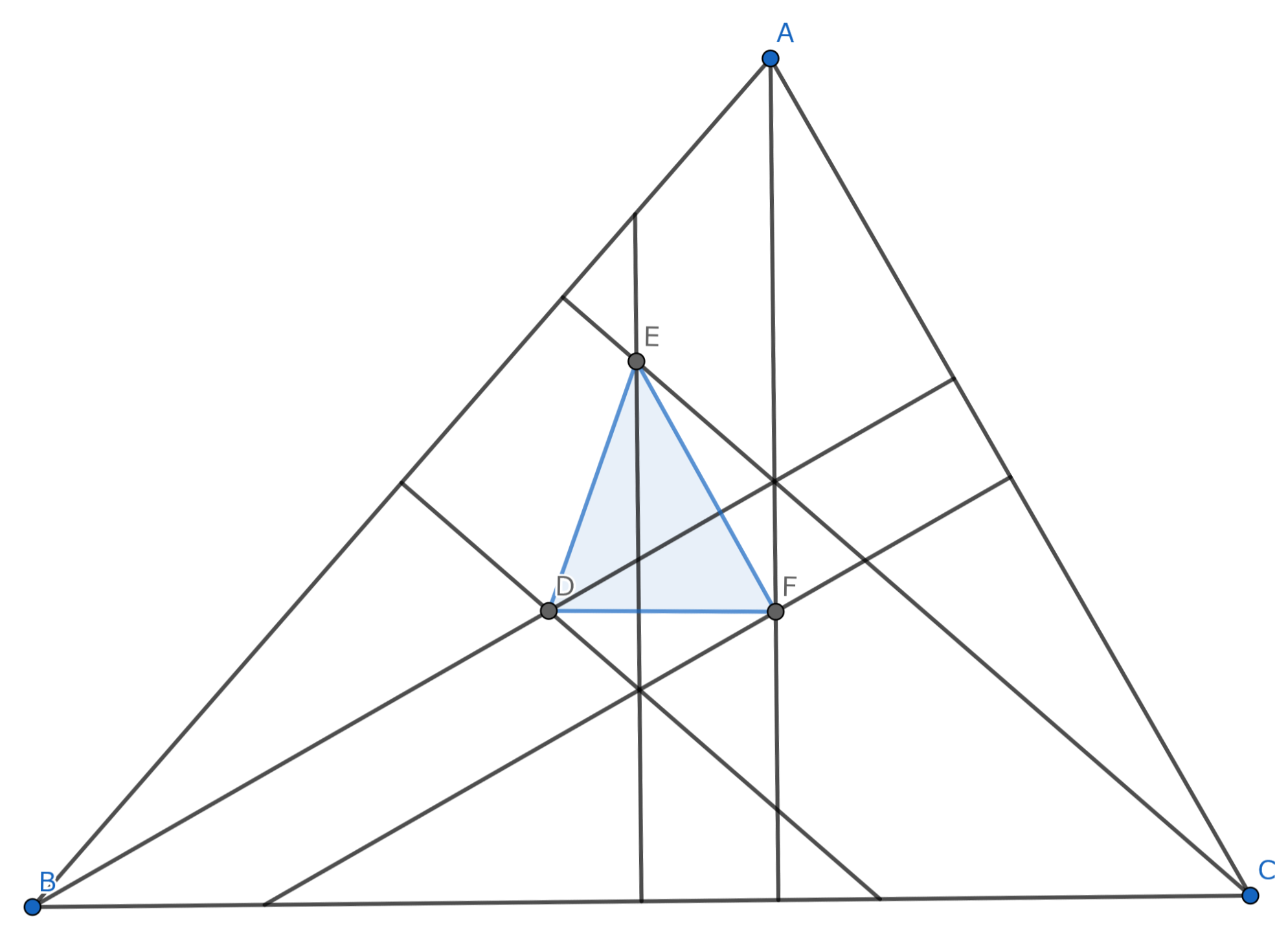}
    \caption{}
    \label{fig:1}
\end{figure}

The following is our main result.
\begin{theorem}\label{MainThm}
    Suppose that \sankaku{ABC} is not an equilateral triangle.
    Let \ten{D} be the intersection of the perpendicular line down from \ten{B} to \ten{AC} and the perpendicular bisector of \ten{AB}, \ten{E} be the intersection of the perpendicular line down from \ten{C} to \ten{AB} and the perpendicular bisector of \ten{BC}, and \ten{F} be the intersection of the perpendicular line down from \ten{A} to \ten{BC} and the perpendicular bisector of \ten{AC} (see Fig. \ref{fig:1}).
    Then, \sankaku{ABC} and \sankaku{DEF} are similar.
\end{theorem}
Note that \sankaku{DEF} degenerates to a single point when \sankaku{ABC} is an equilateral triangle.

We will give the proof of this theorem in two ways.
In the next section, it will be proved in the complex plane.
In section \ref{shotoukika}, we will provide an elementary geometrical proof.

\section{A proof by complex coordinates}

We first introduce an inner product of complex numbers.
Let $i$ be an imaginary unit and $x,y,u,v$ are real numbers.
For complex numbers $\alpha = x+ yi$ and $\beta = u+vi$, we define an \emph{inner product} of $\alpha$ and $\beta$ by
\begin{align}
    \alpha * \beta 
    =xu+uv
    =\frac{\alpha \bar{\beta}+\bar{\alpha}\beta}{2}
\end{align}
where, for a complex number $\zeta$, $\bar{\zeta}$ denotes the complex conjugate of $\zeta$.
Using this notation, for the points $\ten{A,B,C}$ whose complex coordinates are $\alpha ,\beta, \gamma~(\alpha \not= \beta)$ respectively, i.e., $\ten{A}(\alpha), \ten{B}(\beta), \ten{C}(\gamma)$, an equation of a line perpendicular to line \ten{AB} and passing through \ten{C} can be written as
\begin{align}
    (z-\gamma) * (\beta - \alpha) =0
\end{align}
or
\begin{align}
    (z-\gamma)(\bar{\beta}-\bar{\alpha})+(\bar{z}-\bar{\gamma})(\beta - \alpha)=0.
    \label{eq:1}
\end{align}
Then, we can give a proof of Theorem \ref{MainThm}.

\proof[Proof of Theorem \ref{MainThm}]
We fix a complex coordinate system such that points with their own complex coordinates $\ten{O}(0),\ten{A}(\alpha),\ten{B}(\beta),\ten{C}(\gamma)$ satisfy $|\alpha|=|\beta|=|\gamma|=1$.
Note that $\ten{O}$ is the circumcenter of $\sankaku{ABC}$ in this case.
Since the perpendicular bisector of \ten{AB} passes through \ten{O}, its equation can be written from \eqref{eq:1} as
\begin{gather}
    (\beta - \alpha) \bar{z} + (\bar{\beta} - \bar{\alpha})z = 0
    \\
    (\beta - \alpha) \bar{z} + \frac{\alpha - \beta}{\alpha \beta }z = 0
    \\
    z = \alpha \beta \bar{z}.
\end{gather}
The equation of the perpendicular line from \ten{B} down to \ten{AC} can be obtained again from \eqref{eq:1} as
\begin{gather}
    (z-\beta)(\bar{\gamma}-\bar{\alpha})+(\bar{z}-\bar{\beta})(\gamma - \alpha) = 0
    \\
    (z-\beta)\frac{\alpha - \gamma}{\alpha \gamma} + \frac{\beta \bar{z}- 1}{\beta}(\gamma - \alpha) = 0
    \\
    z = \alpha \gamma \bar{z} - \frac{\alpha \gamma}{\beta} + \beta.
\end{gather}
Then, if the coordinate of \ten{D} is denoted by $\delta$, it satisfies
\begin{align}
    \delta = \alpha \beta \bar{\delta}
    \qquad \text{and} \qquad
    \delta = \alpha \gamma \bar{\delta} - \frac{\alpha \gamma }{\beta} + \beta
\end{align}
so that
\begin{gather}
    \beta \delta = \alpha \beta \gamma \bar{\delta} - \alpha \gamma +\beta^2 =\gamma \delta - \alpha \gamma \beta^2 
    \\
    \delta = \frac{\beta^2 - \alpha \gamma}{\beta - \gamma}.
\end{gather}
By the similar arguments, we get
\begin{align}
    D\left( \frac{\beta^2 - \alpha \gamma}{\beta - \gamma} \right),
    \qquad
    E\left( \frac{\gamma^2-\alpha \beta}{\gamma - \alpha} \right),
    \qquad 
    F\left( \frac{\alpha^2 - \beta \gamma}{\alpha - \beta} \right).
\end{align}
Then, the vector $\overrightarrow{\ten{DE}}$ can be represented by
\begin{align}
   \overrightarrow{\ten{DE}} &=  \frac{\gamma^2-\alpha \beta}{\gamma - \alpha} - \frac{\beta^2 - \alpha \gamma}{\beta - \gamma}
   \\
   &= \frac{(\gamma^2 - \alpha \beta )(\beta - \gamma) - (\beta^2 - \alpha \gamma)(\gamma - \alpha)}{(\beta - \gamma) (\gamma - \alpha)}
   \\
   &=\frac{-\gamma (\alpha^2 + \beta^2 + \gamma^2 - \alpha \beta - \beta \gamma - \gamma \alpha)}{(\alpha - \beta)(\beta - \gamma) (\gamma - \alpha)} (\alpha - \beta).
\end{align}
Similarly, we can get
\begin{align}
    \overrightarrow{\ten{EF}} &= \frac{-\alpha (\alpha^2 + \beta^2 + \gamma^2 - \alpha \beta - \beta \gamma - \gamma \alpha)}{(\alpha - \beta)(\beta - \gamma) (\gamma - \alpha)} (\beta - \gamma),
    \\
    \overrightarrow{\ten{FD}} &= \frac{-\beta (\alpha^2 + \beta^2 + \gamma^2 - \alpha \beta - \beta \gamma - \gamma \alpha)}{(\alpha - \beta)(\beta - \gamma) (\gamma - \alpha)} (\gamma - \alpha).
\end{align}
Recalling that $|\alpha|=|\beta|=|\gamma|=1$, we can see
\begin{align}
    \frac{\ten{DE}}{\ten{AB}}=\frac{\ten{EF}}{\ten{BC}}=\frac{\ten{FD}}{\ten{CA}}
    =\left| \frac{\alpha^2 + \beta^2 + \gamma^2 - \alpha \beta - \beta \gamma - \gamma \alpha}{(\alpha - \beta)(\beta - \gamma) (\gamma - \alpha)} \right|
\end{align}
and we can conclude that \sankaku{ABC} and \sankaku{DEF} are similar.
\qed

In the above proof, it has also been shown that the ratio of the similar triangles $\sankaku{ABC}$ and $\sankaku{DEF}$ is
\begin{align}
    \left| \frac{\alpha^2 + \beta^2 + \gamma^2 - \alpha \beta - \beta \gamma - \gamma \alpha}{(\alpha - \beta)(\beta - \gamma) (\gamma - \alpha)} \right|
    =\frac{|\alpha^2 + \beta^2 + \gamma^2 - \alpha \beta - \beta \gamma - \gamma \alpha|}{\ten{AB} \cdot \ten{BC} \cdot \ten{CA}}.
\end{align}

Note that, since a general triangle $\sankaku{XYZ}$ with vertices $(X(\alpha),Y(\beta),Z(\gamma))$ is an equilateral triangle if and only if
\begin{align}
    \alpha^2 + \beta^2 + \gamma^2 - \alpha \beta - \beta \gamma - \gamma \alpha = 0,
\end{align}
it can be seen that \sankaku{DEF} degenerates to a single point if and only if \sankaku{ABC} is an equilateral triangle.

Note also that, if the area of \sankaku{ABC} is denoted by $S$ and the lengths of the sides are $a,b$ and $c$, the ratio of the similar triangles $\sankaku{ABC}$ and $\sankaku{DEF}$ can be represented as 
\begin{align}
    \sqrt{\left( \frac{a^2+b^2+c^2}{8S} \right)^2 - \frac{3}{4}}.
\end{align}

\section{A proof by elementary geometry}\label{shotoukika}

Our proof by elementary geometry is separated by several propositions.

\begin{figure}[t]
    \centering
    \includegraphics[width=0.7\linewidth]{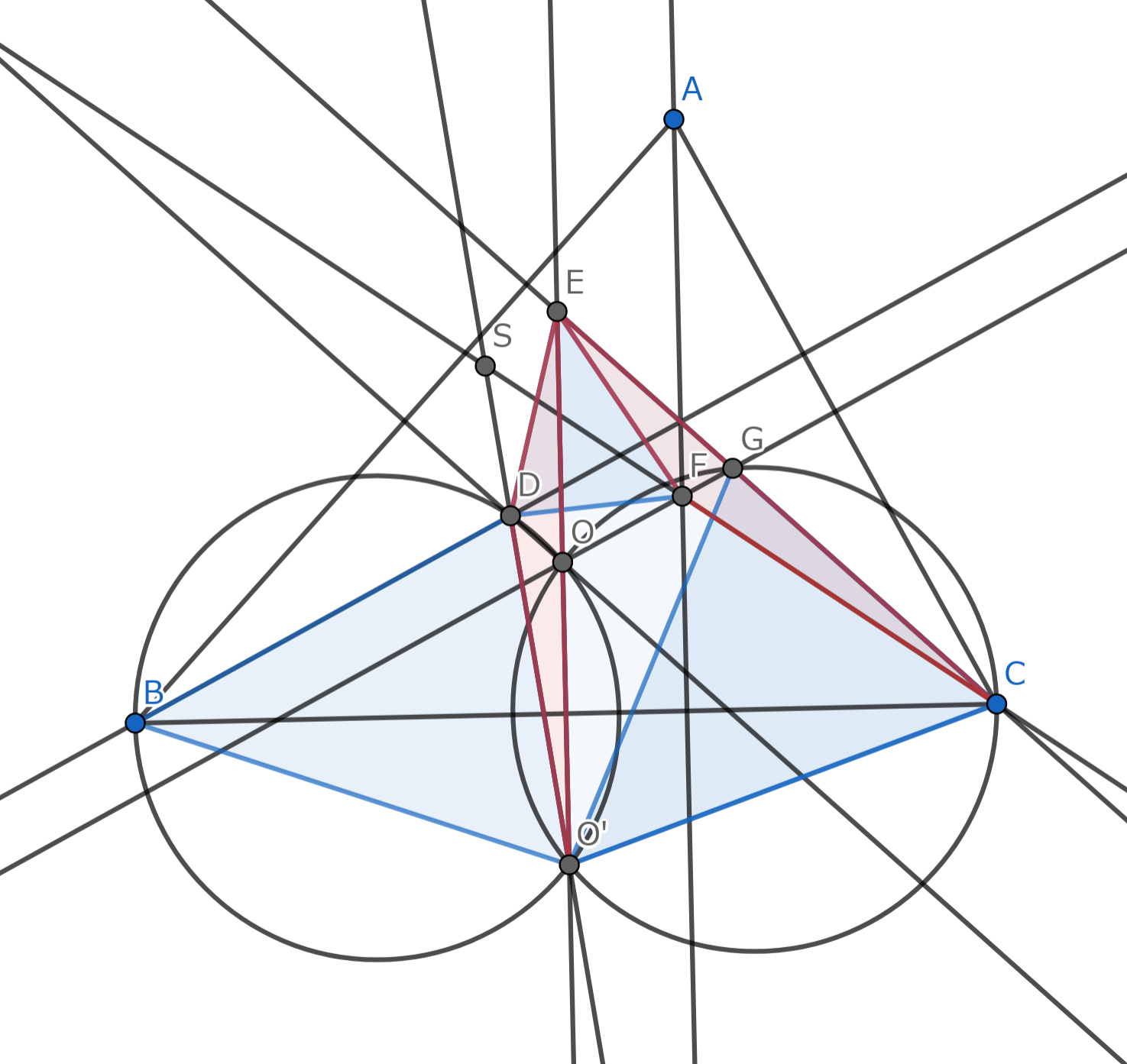}
    \caption{}
    \label{fig:2}
\end{figure}

\begin{proposition}\label{prop1}
    Let \ten{O} be the intersection of the perpendicular bisectors of \ten{AB} and \ten{BC} and \ten{O'} be the point symmetric to \ten{O} with respect to \ten{BC}.
    Then, points \ten{B}, \ten{D}, \ten{O} and \ten{O'} are concyclic.
\end{proposition}
\proof
If a point \ten{P} is taken on the line \ten{OD} that is opposite to \ten{O} with respect to \ten{D}, then
\begin{align}
    \kaku{OO'B}=\kaku{BOO'}=\kaku{BAC}=\kaku{BDP}
\end{align}
holds and the claim is confirmed.
\qed

\begin{proposition}\label{prop2}
    Let \ten{G} be the intersection of the perpendicular line down from \ten{C} to \ten{AB} and the perpendicular bisector of \ten{AC}. 
    Then, points \ten{O}, \ten{G}, \ten{C} and \ten{O'} are concyclic.
\end{proposition}
\proof
Take a point \ten{Q} on the line \ten{OG} on the opposite of \ten{O} with respect to \ten{G}.
Then
\begin{align}
    \kaku{OO'C}=\kaku{COO'}=\kaku{BAC}=\kaku{CGQ}
\end{align}
holds and we can confirm the assertion.
\qed

\begin{proposition}\label{prop3}
    \sankaku{ABC} and \sankaku{O'EC} are similar, where \ten{O'} is the one defined in Proposition \ref{prop1}.
\end{proposition}
\proof
Recognizing that \ten{O} is the circumcenter of \sankaku{ABC},
\begin{align}
    \kaku{CO'E}=\kaku{COO'}=\kaku{BAC}
\end{align}
holds.
Moreover, if \ten{M} denotes the intersection of \ten{OG} and \ten{CA} and \ten{N} is the intersection of \ten{OO'} and \ten{BC}, points \ten{O}, \ten{M}, \ten{C} amd \ten{N} are concyclic.
Then, from Proposition \ref{prop2}, we obtain
\begin{align}
    \kaku{ECO'}=\kaku{EOG}=\kaku{BCA}
\end{align}
and the claim is confirmed.
\qed

\begin{proposition}\label{prop4}
    Let \ten{S} be the intersection of \ten{CF} and \ten{O'D}.
    Then, points \ten{S}, \ten{E}, \ten{C} and \ten{O'} are concyclic.
\end{proposition}
\proof
It suffices to show that $\kaku{ECS}=\kaku{EO'S}$.

From Proposition \ref{prop1}, 
\begin{align}
    \kaku{EO'S}=\kaku{DO'O}=\kaku{DBO}
\end{align}
holds.
We also have
\begin{align}
    \kaku{DBO}=\kaku{DBC}-\kaku{OBC}.
\end{align}
Moreover, it also holds that
\begin{align}
    \kaku{ECS}
    =\kaku{ACS}-\kaku{ACE}
    =\kaku{ACF}-\kaku{ACE}
    =\kaku{CAF}-\kaku{ACE}
\end{align}
by the definition of \ten{S} and 
\begin{align}
    \kaku{OBC}=\kaku{ACE},
    \qquad 
    \kaku{DBC}=\kaku{CAF}.
\end{align}
From these facts, we can conclude that $\kaku{ECS}=\kaku{EO'S}$.
\qed

\begin{proposition}\label{prop5}
\sankaku{DBO'} and \sankaku{GO'C} are similar to \sankaku{ABC}.
Furthermore, \sankaku{DBO'} and \sankaku{GO'C} are congruent.
\end{proposition}
\proof
The former can be derived immediately from Proposition \ref{prop1} and \ref{prop2} and a little angle calculation.
The latter is clear from $\ten{BO'}=\ten{O'C}$.
\qed

\begin{proposition}\label{prop6}
    \sankaku{ECF} and \sankaku{EO'D} are similar.
\end{proposition}
\proof
Since it is shown in Proposition \ref{prop4} that $\kaku{ECS}=\kaku{EO'S}$, it suffices to show 
\begin{align}
    \frac{\ten{EC}}{\ten{EO'}} = \frac{\ten{CF}}{\ten{O'D}}.
\end{align}
From Proposition \ref{prop3}, we have
\begin{align}
    \frac{\ten{EC}}{\ten{EO'}} = \frac{\ten{BC}}{\ten{BA}}
\end{align}
and, by Proposition \ref{prop5},
\begin{align}
    \frac{\ten{CF}}{\ten{O'D}}=\frac{\ten{CF}}{\ten{CG}}=\frac{\ten{AF}}{\ten{CG}}.
\end{align}
Moreover, 
\begin{align}
    \frac{\ten{AF}}{\ten{CG}}=\frac{2 \ten{AF}}{2 \ten{CG}}=\frac{\ten{AC}/\sin \kaku{BCA}}{\ten{AC}/\sin\kaku{CAB}}=\frac{\sin \kaku{CAB}}{\sin \kaku{BCA}}
\end{align}
and, by the law of sines,
\begin{align}
    \frac{\sin \kaku{CAB}}{\sin \kaku{BCA}} = \frac{\ten{BC}}{\ten{BA}}.
\end{align}
Then, we get the assertion.
\qed

\proof[Proof of Theorem \ref{MainThm}]

From Proposition \ref{prop6}, we see that \sankaku{DEF} and \sankaku{O'EC} are similar (by spiral similarity) and \sankaku{O'EC} and \sankaku{ABC} are similar by Proposition \ref{prop3}, which finishes the proof.
\qed

\section{Further extensions}

The theorem proved in the previous section states that the intersections of the perpendicular bisectors of each side of the triangle and the perpendicular line from each vertex to its opposite side create a triangle similar to the original one.
In this section, we show that, by considering the intersections of the perpendicular bisectors of each side and lines perpendicular to sides of the triangle passing through one of the vertices, a new similar triangle emerges.

\begin{figure}[t]
    \centering
    \includegraphics[width=0.7\linewidth]{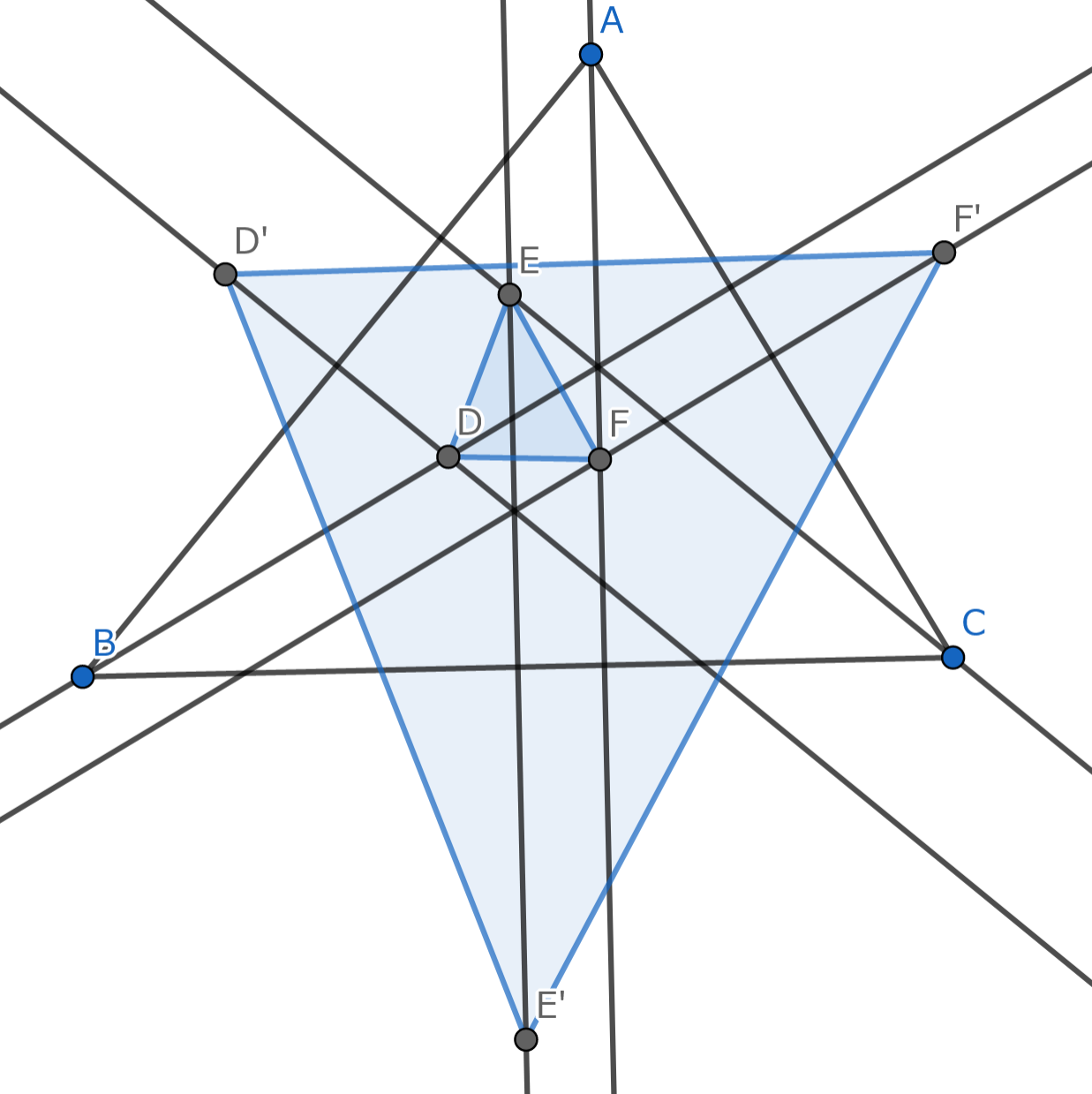}
    \caption{}
    \label{fig:7left}
\end{figure}

Let \ten{D'} be the intersection of the perpendicular bisector of \ten{AB} and the line passing through \ten{A} and perpendicular to \ten{CA} (see Fig. \ref{fig:7left}).
Then, the complex coordinate of \ten{D'} is
\begin{align}
    D'\left( \beta \frac{\alpha - \gamma}{\beta - \gamma} \right).
\end{align}
Note that \ten{D'} is also a point symmetric to \ten{D} with respect to \ten{AB}.
Similarly, \ten{E'} being the intersection of the perpendicular bisector of \ten{BC} and the line passing through \ten{B} and perpendicular to \ten{AB} and \ten{F'} being the intersection of the perpendicular bisector of \ten{CA} and the line passing through \ten{C} and perpendicular to \ten{BC}, the complex coordinates of \ten{E'} and \ten{F'} are
\begin{align}
    E'\left( \gamma \frac{\beta - \alpha}{\gamma - \alpha} \right),
    \qquad 
    F'\left( \alpha \frac{\gamma - \beta }{\alpha - \beta} \right).
\end{align}
Then, since we get
\begin{align}
    \overrightarrow{\ten{D'E'}} = \gamma \frac{\beta - \alpha}{\gamma - \alpha} - \beta \frac{\alpha - \gamma}{\beta - \gamma} = \frac{\alpha^2 \beta + \beta^2 \gamma + \gamma^2 \alpha - 3\alpha \beta \gamma}{(\alpha - \beta )(\beta - \gamma)(\gamma - \alpha)}(\alpha - \beta)
\end{align}
and
\begin{align}
    \overrightarrow{\ten{E'F'}} &= \frac{\alpha^2 \beta + \beta^2 \gamma + \gamma^2 \alpha - 3\alpha \beta \gamma}{(\alpha - \beta )(\beta - \gamma)(\gamma - \alpha)}(\beta - \gamma),
    \\
    \overrightarrow{\ten{F'D'}} &= \frac{\alpha^2 \beta + \beta^2 \gamma + \gamma^2 \alpha - 3\alpha \beta \gamma}{(\alpha - \beta )(\beta - \gamma)(\gamma - \alpha)}(\gamma - \alpha),
\end{align}
in the similar way, we see that \sankaku{D'E'F'} and \sankaku{ABC} are similar.
Noting that their similarity ratio is
\begin{align}
    \left| \frac{\alpha^2 \beta + \beta^2 \gamma + \gamma^2 \alpha - 3\alpha \beta \gamma}{(\alpha - \beta )(\beta - \gamma)(\gamma - \alpha)} \right|,
\end{align}
since 
\begin{align}
    &|\alpha^2 \beta + \beta^2 \gamma + \gamma^2 \alpha - 3\alpha \beta \gamma|^2
    \\
    &=|(\alpha - \beta )(\beta - \gamma)(\gamma - \alpha)|^2 + |\alpha^2 + \beta^2 + \gamma^2 - \alpha\beta - \beta\gamma - \gamma\alpha|^2
\end{align}
holds, we obtain
\begin{align}
    \left| \frac{\alpha^2 \beta + \beta^2 \gamma + \gamma^2 \alpha - 3\alpha \beta \gamma}{(\alpha - \beta )(\beta - \gamma)(\gamma - \alpha)} \right|^2 
    = 1 + \left| \frac{\alpha^2 + \beta^2 + \gamma^2 - \alpha\beta - \beta\gamma - \gamma\alpha}{(\alpha - \beta )(\beta - \gamma)(\gamma - \alpha)} \right|^2.
\end{align}
This shows that the area of \sankaku{D'E'F'} is equal to the sum of the area of \sankaku{ABC} and \sankaku{DEF}.

Please note that 
\begin{align}
    \frac{\alpha^2 \beta + \beta^2 \gamma + \gamma^2 \alpha - 3\alpha \beta \gamma}{(\alpha - \beta )(\beta - \gamma)(\gamma - \alpha)}
    =\frac{\alpha\beta\gamma \left( \dfrac{\alpha}{\gamma} + \dfrac{\beta}{\alpha} + \dfrac{\gamma}{\beta} - 3 \right)}{(\alpha - \beta )(\beta - \gamma)(\gamma - \alpha)}
     =\frac{\alpha\beta\gamma \overline{\dfrac{\gamma}{\alpha} + \dfrac{\alpha}{\beta} + \dfrac{\beta}{\gamma} - 3}}{(\alpha - \beta )(\beta - \gamma)(\gamma - \alpha)}
\end{align}
is a cyclic formula but not a symmetric formula.
However, if we take the absolute value the above formula, it becomes
\begin{align}
    \left| \frac{\alpha\beta\gamma \overline{\dfrac{\gamma}{\alpha} + \dfrac{\alpha}{\beta} + \dfrac{\beta}{\gamma} - 3}}{(\alpha - \beta )(\beta - \gamma)(\gamma - \alpha)} \right|
    = \left| \frac{\alpha\beta\gamma \left(\dfrac{\gamma}{\alpha} + \dfrac{\alpha}{\beta} + \dfrac{\beta}{\gamma} - 3\right)}{(\alpha - \beta )(\beta - \gamma)(\gamma - \alpha)} \right|
\end{align}
so that the value will not change no matter which two variables are switched.

\begin{figure}[t]
      \centering
    \includegraphics[width=0.7\linewidth]{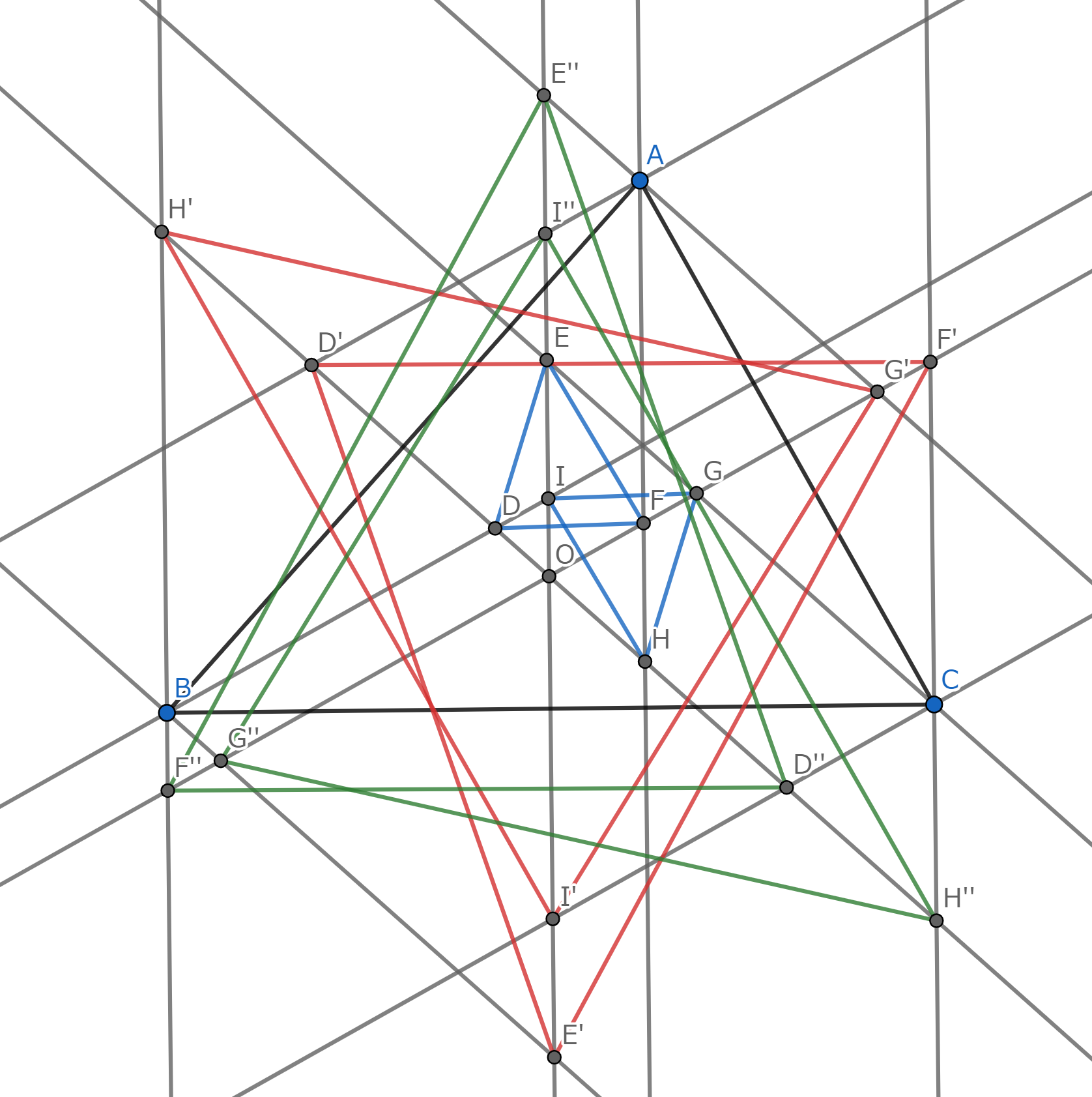}
    \caption{}
    \label{fig:7right}
\end{figure}

Note also that, as seen in Fig. \ref{fig:7right}, by properly choosing the intersections of perpendicular lines and perpendicular bisectors on the triangles, we can create four triangles congruent to \sankaku{D'E'F'}. 

\begin{thebibliography}{9}


\bibitem{Naka}
Naka, H.
``A theorem found using a drawing app (in Japanese)."
S\=ugaku seminar. \textbf{59} (11) (2020), p. 56, NIPPON HYORON SHA CO., Tokyo. 

\bibitem{Naka2022}
Naka, H.
``A theorem in plane geometry and its relation to Miquel points (in Japanese)."
Naomichi Shiono Memorial 10th ``Independent Research in Arithmetic and Mathematics" Contest 2022 Winning Works
\url{https://www.rimse.or.jp/research/past/winner10th.html}.
Available at:
\url{https://www.rimse.or.jp/research/past/pdf/10th/work06.pdf}

\bibitem{Miquel}
Weisstein, E. W.
``Miquel's Theorem." From MathWorld--A Wolfram Web Resource. \url{https://mathworld.wolfram.com/MiquelsTheorem.html}



\end{thebibliography}
\end{document}